\def\draft{n}
\documentclass{amsart}
\usepackage[headings]{fullpage}
\usepackage{amssymb,epic,eepic,epsfig}
\usepackage{graphicx}
\usepackage{texdraw}
\usepackage{url}
\usepackage[bookmarks=true,%
    colorlinks=true,%
    linkcolor=blue,%
    citecolor=blue,%
    filecolor=blue,%
    menucolor=blue,%
    urlcolor=blue,%
    breaklinks=true]{hyperref}

\newtheorem{theorem}{Theorem}[section]
\newtheorem{proposition}{Proposition}[section]
\theoremstyle{definition}

\newtheorem{definition}[proposition]{Definition}
\newtheorem{remark}[proposition]{Remark}

\def\printname#1{
        \if\draft y
                \smash{\makebox[0pt]{\hspace{-0.5in}
                        \raisebox{8pt}{\tt\tiny #1}}}
        \fi
}

\newcommand{\psdraw}[2]
         {\begin{array}{c} \hspace{-1.3mm}
        \raisebox{-4pt}{\epsfig{figure=draws/#1.eps,width=#2}}
        \hspace{-1.9mm}\end{array}}

\newlength{\standardunitlength}
\catcode`\@=11
\long\def\@makecaption#1#2{%
     \vskip 10pt

\setbox\@tempboxa\hbox{
       \small\sf{\bfcaptionfont #1. }\ignorespaces #2}%
     \ifdim \wd\@tempboxa >\captionwidth {%
         \rightskip=\@captionmargin\leftskip=\@captionmargin
         \unhbox\@tempboxa\par}%
       \else
         \hbox to\hsize{\hfil\box\@tempboxa\hfil}%
     \fi}
\font\bfcaptionfont=cmssbx10 scaled \magstephalf
\newdimen\@captionmargin\@captionmargin=2\parindent
\newdimen\captionwidth\captionwidth=\hsize
\catcode`\@=12

\def\lbl#1{\label{#1}\printname{#1}}

\def\BN{\mathbb N}
\def\BZ{\mathbb Z}

\def\BQ{\mathbb Q}

\def\BC{\mathbb C}

\def\l{\lambda}

\def\ga{\gamma}

\def\e{\epsilon}

\def\longto{\longrightarrow}

\def\SL{\mathrm{SL}}

\def\GL{\mathrm{GL}}
\def\NP{\mathrm{NP}}
\def\NQ{\mathrm{NQ}}

\begin{document}


\title[The $A$-polynomial of the $(-2,3,3+2n)$ pretzel knots]{
The $A$-polynomial of the $(-2,3,3+2n)$ pretzel knots}
\author{Stavros Garoufalidis}
\address{School of Mathematics \\
         Georgia Institute of Technology \\
         Atlanta, GA 30332-0160, USA \newline
         {\tt \url{http://www.math.gatech.edu/~stavros }}}
\email{stavros@math.gatech.edu}
\author{Thomas W. Mattman}
\address{Department of Mathematics and Statistics \\
        California State University, Chico \\
        Chico, CA 95929-0525, USA \newline
        {\tt  \url{http://www.csuchico.edu/~tmattman/}}}
\email{TMattman@CSUChico.edu}
\thanks{S.G. was supported in part by NSF. \\
\newline
1991 {\em Mathematics Classification.} Primary 57N10. Secondary 57M25.
}

\date{April 8, 2011}


\begin{abstract}
We show that the $A$-polynomial $A_n$ of the 1-parameter family of pretzel 
knots $K_n=(-2,3,3+2n)$ satisfies a linear recursion relation of order 4 with
explicit constant coefficients and initial conditions. Our proof combines
results of Tamura-Yokota and the second author. As a corollary, we show
that the $A$-polynomial of $K_n$ and the mirror of $K_{-n}$ are related
by an explicit $\GL(2,\BZ)$ 
action. We leave open the question of whether or not 
this action lifts to the quantum level.
\end{abstract}

\maketitle


\section{Introduction}
\lbl{sec.intro}

\subsection{The behavior of the $A$-polynomial under filling}
\lbl{sub.filling}

In \cite{CCGLS}, the authors introduced the $A$-polynomial $A_W$ of a 
hyperbolic 3-manifold $W$ with one cusp.
It is a 2-variable polynomial which describes
the dependence of the eigenvalues of a meridian and longitude
under any representation of $\pi_1(W)$ into $\SL(2,\BC)$. The $A$-polynomial 
plays a key role in two problems:
\begin{itemize}
\item
the deformation of the hyperbolic structure of $W$,
\item
the problem of exceptional (i.e., non-hyperbolic) fillings of $W$.
\end{itemize}
Knowledge of the $A$-polynomial (and often, of its Newton polygon)
is translated directly into information about the above problems, and 
vice-versa.
In particular, as demonstrated by Boyer and Zhang \cite{BZ}, the Newton 
polygon is dual
to the fundamental polygon of the Culler-Shalen seminorm \cite{CGLS} and, 
therefore, can be used to classify
cyclic and finite exceptional surgeries.

In \cite{Ga3}, the first author observed a pattern in the behavior 
of the $A$-polynomial (and its Newton polygon) of a 1-parameter family of 
3-manifolds obtained by fillings of a 2-cusped manifold. To state the
pattern, we need to introduce some notation. 
Let $K=\BQ(x_1,\dots,x_r)$ denote the field of rational functions in $r$ 
variables $x_1,\dots,x_r$.

\begin{definition}
\lbl{def.holo}
We say that a sequence of rational functions $R_n \in K$ 
(defined for all integers $n$) is {\em holonomic} if it satisfies a linear 
recursion with constant coefficients. In other words, there exists a natural
number $d$ and $c_k \in K$ for $k=0,\dots,d$ with $c_d c_0 \neq 0$
such that for all integers $n$ we have:
\begin{equation}
\lbl{eq.recRR}
\sum_{k=0}^d c_k R_{n+k}=0
\end{equation}
\end{definition}
Depending on the circumstances, one can restrict attention to sequences
indexed by the natural numbers (rather than the integers). 

Consider a 
hyperbolic manifold $W$ with two cusps $C_1$ and $C_2$. 
Let $(\mu_i,\l_i)$ for $i=1,2$ be pairs of 
meridian-longitude curves, and let $W_n$ denote the result of $-1/n$ filling 
on $C_2$. Let $A_n(M_1,L_1)$ denote the $A$-polynomial of $W_n$ with the
meridian-longitude pair inherited from $W$.

\begin{theorem}\cite{Ga3}
\lbl{thm.Ga1}
With the above conventions, there exists a holonomic sequence 
$R_n(M_1,L_1) \in \BQ(M_1,L_1)$ such that 
for all but finitely many integers $n$, $A_n(M_1,L_1)$ divides the 
numerator of $R_n(M_1,L_1)$.
In addition, a recursion for $R_n$ can be computed explicitly via elimination,
from an ideal triangulation of $W$.
\end{theorem}

\subsection{The Newton polytope of a holonomic sequence}
\lbl{sub.newton}

Theorem \ref{thm.Ga1} motivates us to study the Newton polytope of a 
holonomic sequence of Laurent polynomials. To state our result,
we need some definitions. Recall that the {\em Newton polytope} of
a Laurent polynomial in $n$ variables $x_1,\dots,x_n$ is the convex hull
of the points whose coordinates are the exponents of its monomials. 
Recall that a {\em quasi-polynomial} is a function $p:\BN\longto\BQ$ of
the form $p(n)=\sum_{k=0}^d c_k(n) n^k$ where $c_k: \BN \longto \BQ$ are
periodic functions. When $c_d \neq 0$, we call $d$ the {\em degree}
of $p(n)$. We will call quasi-polynomials of degree at most one (resp. 
two) {\em quasi-linear} (resp. {\em quasi-quadratic}).
Quasi-polynomials
appear in lattice point counting problems (see \cite{Eh,CW}), in the 
Slope Conjecture in quantum topology (see \cite{Ga1}), in enumerative
combinatorics (see \cite{Ga2}) and also in the $A$-polynomial of filling
families of 3-manifolds (see \cite{Ga3}).

\begin{definition}
We say that a sequence $N_n$ of polytopes is linear 
(resp. quasi-linear) if the coordinates of the vertices of $N_n$ are
polynomials (resp. quasi-polynomials) in $n$ of degree at most one. Likewise,
we say that a sequence  $N_n$ of polytopes is quadratic 
(resp. quasi-quadratic) if the coordinates of the vertices of $N_n$ are
polynomials (resp. quasi-polynomials) of degree at most two.
\end{definition}

\begin{theorem}\cite{Ga3}
\lbl{thm.Ga2}
Let $N_n$ be the Newton polytope of a holonomic sequence $R_n \in
\BQ[x_1^{\pm 1},\dots,x_r^{\pm 1}]$. Then, for all but finitely many integers $n$,
$N_n$ is quasi-linear.
\end{theorem}

\subsection{Do favorable links exist?}
\lbl{sub.favor}

Theorems \ref{thm.Ga1} and \ref{thm.Ga2} are general, but in favorable
circumstances more is true. Namely, consider a family of knot complements 
$K_n$, obtained by $-1/n$ filling on a cusp of a 2-component hyperbolic 
link $J$. Let $f$ denote the linking number of the two components of $J$, 
and let $A_n$ denote the $A$-polynomial of $K_n$ with respect to
its canonical meridian and longitude $(M,L)$. By definition, $A_n$
contains all components of irreducible representations, but {\em not}
the component $L-1$ of abelian representations.

\begin{definition}
\lbl{def.favor}
We say that $J$, a 2-component link in 3-space, with
linking number $f$ is {\em favorable} if 
$A_n(M,LM^{-f^2 n}) \in \BQ[M^{\pm 1},L^{\pm 1}]$ is holonomic.
\end{definition}
The shift of coordinates, $L M^{-f^2 n}$, above 
is due to the canonical meridian-longitude pair of $K_n$
differing from the corresponding pair for the unfilled component of $J$
as a result of the nonzero linking number. Theorem \ref{thm.Ga2} combined 
with the above shift implies that, for a favorable link, 
the Newton polygon of $K_n$ is quasi-quadratic.

Hoste-Shanahan studied the first examples of a favorable link, the
{\em Whitehead link} and its {\em half-twisted version} 
(see Figure \ref{3links}), and consequently 
gave an explicit recursion relation for the 1-parameter families of 
$A$-polynomials of twist knots $K_{2,n}$ and $K_{3,n}$ respectively; 
see \cite{HS}. 

\begin{figure}[htpb]
$$ 
\psdraw{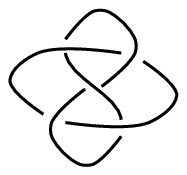}{1.5in} \qquad
\psdraw{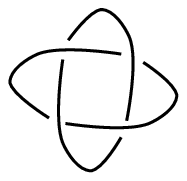}{1.5in} \qquad
\psdraw{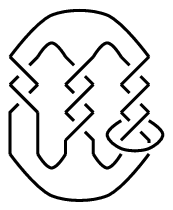}{1.3in}
$$
\caption{The Whitehead link on the left, the half-twisted Whitehead link
in the middle and our seed link $J$ at right.}\lbl{3links}
\end{figure}

The goal of our paper is to give another example of a favorable link $J$ 
(see Figure \ref{3links}),
whose 1-parameter filling gives rise to the family of $(-2,3,3+2n)$
{\em pretzel knots}. Our paper is a concrete illustration of the general
Theorems \ref{thm.Ga1} and \ref{thm.Ga2} above.  Aside from this, the
1-parameter family of knots $K_n$, where $K_n$ is the $(-2,3,3+2n)$ pretzel 
knot,
is well-studied in hyperbolic geometry (where $K_n$ and the mirror of 
$K_{-n}$ are pairs of geometrically
similar knots; see \cite{BH,MM}), in exceptional Dehn surgery (where
for instance $K_2=(-2,3,7)$ has three Lens space fillings $1/0$, $18/1$
and $19/1$; see \cite{CGLS}) and in Quantum Topology (where
$K_n$ and the mirror of $K_{-n}$ have different Kashaev invariant, equal
volume, and different subleading corrections to the volume, see
\cite{GZ}).
 
The success of Theorems \ref{thm.1} and \ref{thm.2} below hinges on two 
independent results of Tamura-Yokota and the second author \cite{TY,Ma},
and an additional lucky coincidence.
Tamura-Yokota compute an explicit recursion relation, as in Theorem
\ref{thm.1}, by elimination, using the gluing equations of the decomposition 
of the complement of $J$ into six ideal tetrahedra; see \cite{TY}.
The second author computes the Newton polygon $N_n$ of the $A$-polynomial 
of the family $K_n$ of pretzel knots; see \cite{Ma}. 
This part is considerably more difficult, and requires: 
\begin{itemize}
\item[(a)] The set of boundary slopes of $K_n$, which are available by
applying the Hatcher-Oertel algorithm \cite{HO,D} to the 1-parameter family 
$K_n$
of Montesinos knots. The four slopes given by the algorithm are candidates 
for the slopes of the sides of $N_n$. 
Similarly,  the fundamental polygon of the
Culler-Shalen seminorm of $K_n$ has vertices in rays
which are the multiples of the slopes of $N_n$. Taking advantage of the 
duality of the fundamental polygon and Newton polygon, in order to describe 
$N_n$ it is enough to determine
the vertices of the Culler-Shalen polygon. 
\item[(b)] Use of the exceptional $1/0$ filling and 
two fortunate exceptional Seifert fillings of $K_n$ with slopes $4n+10$ and 
$4n+11$ to determine exactly the vertices of the 
Culler-Shalen polygon and consequently $N_n$.
In particular, the boundary slope $0$ is not a side of $N_n$ 
(unless $n = -3$) and
the Newton polygon is a hexagon for all hyperbolic $K_n$. 
\end{itemize}
Given the work of \cite{TY} and \cite{Ma}, 
if one is lucky enough to match $N_n$ of \cite{Ma} with the Newton polygon
of the solution of the recursion relation of \cite{TY} (and also match
a leading coefficient), then Theorem \ref{thm.1} below follows;
i.e., $J$ is a favorable link.

\subsection{Our results for the pretzel knots $K_n$}
\lbl{sub.results}

Let $A_n(M,L)$ denote the $A$-polynomial of the pretzel knot $K_n$, using
the canonical meridian-longitude coordinates. 
Consider the sequences of Laurent polynomials $P_n(M,L)$ and $Q_n(M,L)$
defined by:
\begin{equation}
\lbl{eq.Pn}
P_n(M,L)=A_n(M,LM^{-4n}) 
\end{equation}
for $n>1$ and 
\begin{equation}
\lbl{eq.Qn}
Q_n(M,L)=A_n(M,LM^{-4n}) M^{-4(3n^2+11n+4)}
\end{equation}
for $n<-2$ and $Q_{-2}(M,L)=A_{-2}(M,LM^{-8}) M^{-20}$. 
In the remaining cases 
$n=-1,0,1$, the knot $K_n$ is not hyperbolic (it is the torus knot
$5_1$, $8_{19}$ and $10_{124}$ respectively), and one expects exceptional 
behavior. This is reflected in the fact that $P_n$ for $n=0,1$ and
$Q_n$ for $n=-1,0$ can be defined to be suitable rational functions 
(rather than polynomials) of $M,L$.
Let $\NP_n$ and $\NQ_n$ denote the Newton polygons of $P_n$ and $Q_n$
respectively. 

\begin{theorem}
\lbl{thm.1}
\rm{(a)}
$P_n$ and $Q_n$ satisfy linear recursion relations
\begin{equation}
\lbl{eq.recP}
\sum_{k=0}^4 c_k P_{n+k}=0, \qquad n \geq 0 
\end{equation}
and
\begin{equation}
\lbl{eq.recQ}
\sum_{k=0}^4 c_{k} Q_{n-k}=0,  \qquad n \leq 0 
\end{equation}
where the coefficients $c_k$ and the initial conditions $P_n$ for
$n=0,\dots,3$ and $Q_n$ for $n=-3,\dots,0$ are given in Appendix \ref{sec.ck}.
\newline
\rm{(b)} In $(L,M)$ coordinates, $\NP_n$ and $\NQ_n$ are hexagons
with vertices 
\begin{equation}
\lbl{eq.hexagon1}
\{ \{0, 0\}, \{1, - 4 n+16\}, \{ n-1,  12n-12 \}, \{2 n+1, 
 16 n+18 \}, \{3 n-1, 32 n-10\}, \{3 n, 28 n+6 \}\}
\end{equation}
for $P_n$ with $n>1$ 
\begin{figure}[htpb]
$$ 
\psdraw{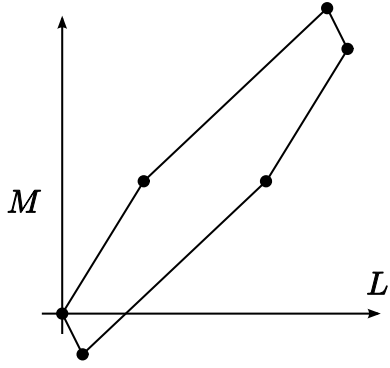}{2in} 
$$
\caption{The Newton polygon $\NP_n$.}
\end{figure}
and 
\begin{equation}
\lbl{eq.hexagon2}
%
%
\{\{0, 4 n+28 \}, \{1, 38\},
 \{-n,  - 12 n +26 \},
 \{- 2 n - 3, -16 n-4 \},
 \{-3 n - 4, -28 n -16 \}, \{-3 n -3, - 32 n-6 \}\}
\end{equation}
for $Q_n$ with $n<-1$.
\end{theorem}

\begin{remark}
\lbl{rem.1}
We can give a single recursion relation valid for 
$n \in \BZ\setminus\{-1,0,1\}$ as follows. Define
\begin{equation}
\lbl{eq.Rn}
R_n(M,L)=A_n(M,LM^{-4n}) b^{|n|} \e_n(M),
\end{equation}
where
\begin{equation}
\lbl{eq.b}
b=\frac{1}{L M^8 (1-M^2)(1+ L M^{10})}
\qquad c=\frac{L^3 M^{12} (1 - M^2)^3}{(1 + L M^{10})^3} 
\quad \e_n(M)=
\begin{cases} 
1 & \text{if $n> 1$} \\
c M^{-4 (3 + n) (2 + 3 n)} & \text{if $n<-2$} \\
c M^{-28} & \text{if $n=-2$}
\end{cases}
\end{equation}
Then, $R_n$ satisfies the palindromic fourth order linear recursion
\begin{equation}
\lbl{eq.recR}
\sum_{k=0}^4 \ga_k R_{n+k}=0
\end{equation}
where the coefficients $\ga_k$ and the initial conditions $R_n$
for $n=0,\dots,3$ are given in Appendix \ref{sec.Rn}. Moreover,
$R_n$ is related to $P_n$ and $Q_n$ by:
\begin{equation}
\lbl{eq.PQR}
R_n = \begin{cases}
P_n b^{|n|} & \text{if $n \geq 0$} \\
Q_n b^{|n|}c M^{-8} & \text{if $n \leq 0$} 
\end{cases}
\end{equation}
\end{remark}

\begin{remark}
\lbl{rem.3}
The computation of the Culler-Shalen seminorm of the pretzel knots $K_n$
has an additional application, namely it determines the number
of components (containing the character of an irreducible representation)
of the $\SL(2,\BC)$ character variety of the knot, and 
consequently the number of factors of its $A$-polynomial. In the case
of $K_n$, (after translating the results of \cite{Ma} for the 
pretzel knots $(-2,3,n)$ to the pretzel knots $(-2,3,3+2n)$)
it was shown by the second author 
\cite[Theorem 1.6]{Ma} that the character variety of $K_n$ has 
one (resp. two) components when $3$ does not divide $n$ (resp. divides $n$).
The non-geometric factor of $A_n$ is given by 
$$
\begin{cases}
1-L M^{4(n+3)} & n \geq 3 \\
L-M^{-4(n+3)} & n \leq -3
\end{cases}
$$
for $n \neq 0$ a multiple of $3$.
\end{remark}

Since the $A$-polynomial has even powers of $M$, we can define
the $B$-polynomial by
$$
B(M^2,L)=A(M,L).
$$
Our next result relates the $A$-polynomials of the geometrically similar pair
$(K_n,-K_{-n})$ by an explicit $\GL(2,\BZ)$ transformation.

\begin{theorem}
\lbl{thm.2}
For $n>1$ we have:
\begin{equation}
\lbl{eq.duality}
B_{-n}(M,LM^{2n-5}) 
= (-L)^n M^{ 3(2n^2-7n+7)} B_n(-L^{-1}, L^{2n+5} M^{-1}) \eta_n
\end{equation}
where $\eta_n=1$ (resp. $M^{22}$) when $n>2$ (resp. $n=2$).
\end{theorem}

\section{Proofs}
\lbl{sec.proofs}

\subsection{The equivalence of Theorem \ref{thm.1} and Remark \ref{rem.1}}
\lbl{sec.rem1}

In this subsection we will show the equivalence of Theorem \ref{thm.1}
and Remark \ref{rem.1}. Let $\ga_k=c_k/b^k$ for $k=0,\dots,4$
where $b$ is given by \eqref{eq.b}. It is easy to see that the $\ga_k$
are given explicitly by Appendix \ref{sec.Rn}, and moreover, they
are palindromic. Since $R_n=P_n b^n$ for $n=0,\dots,3$ it follows
that $R_n$ and $P_n b^n$ satisfy the same recursion relation \eqref{eq.recR}
for $n \geq 0$ with the same initial conditions. It follows that
$R_n=P_n b^n$ for $n \geq 0$.

Solving \eqref{eq.recR} backwards, we can check by an explicit calculation
that $R_n = Q_n b^{|n|} c M^{-8}$ for $n=-3,\dots,0$ where $b$ and $c$ are given 
by \eqref{eq.b}. Moreover, $R_n$ and $Q_n b^{|n|} c M^{-8}$ satisfy the same 
recursion relation \eqref{eq.recR} for $n<0$. It follows that 
$R_n = Q_n b^{|n|} c M^{-8}$ for $n<0$.
This concludes the proof of Equations \eqref{eq.recR} and \eqref{eq.PQR}.

\subsection{Proof of Theorem \ref{thm.1}}
\lbl{sub.thm1}

Let us consider first the case of $n \geq 0$, and denote by
$P'_n$ for $n \geq 0$ the unique solution to the
linear recursion relation \eqref{eq.recP} with the
initial conditions as in Theorem \ref{thm.1}. Let $R'_n=P'_n b^n$ be 
defined according to Equation \eqref{eq.PQR} for $n \geq 0$.

Remark \ref{rem.1} implies that $R'_n$ satisfies the recursion relation 
of \cite[Thm.1]{TY}. It follows by \cite[Thm.1]{TY} that $A_n(M,LM^{-4n})$ 
divides $P'_n(M,L)$ when $n>1$.

Next, we claim that the Newton polygon $\NP'_n$ of $P'_n(M,L)$ is given by 
\eqref{eq.hexagon1}. This can be verified easily by induction on $n$. 

Next, in \cite[p.1286]{Ma}, the second author computes the Newton polygon 
$N_n$ of the $A_n(M,L)$. It is a hexagon given in $(L,M)$ coordinates
by
\begin{equation*}
\begin{split}
\{ & \{0, 0\}, \{1, 16\}, \{n - 1, 4 (n^2 + 2 n - 3)\}, \{2 n + 1, 
  2 (4 n^2 + 10 n + 9)\}, \\ & \{3 n - 1, 2 (6 n^2 + 14 n - 5)\}, \{3 n, 
  2 (6 n^2 + 14 n + 3)\}\}
\end{split}
\end{equation*}
when $n>1$,
\begin{equation*}
\begin{split}
\{ & \{-3 n - 4, 0\}, \{-3 (1 + n), 10\}, \{-3 - 2 n, 4 (3 + 4 n + n^2)\}, \\
& \{-n,  2 (4 n^2 + 16 n + 21)\}, \{0, 4 (3 n^2 + 12 n  + 11)\}, \{1, 
  6 (2 n^2 + 8 n + 9)\}\}
\end{split}
\end{equation*}
when $n<-2$
and
$$
\{\{0, 0\}, \{1, 0\}, \{2, 4\}, \{1, 10\}, \{2, 14\}, \{3, 14\}\}
$$
when $n=-2$.
Notice that the above 1-parameter families of Newton polygons are quadratic.
It follows by explicit calculation that the Newton polygon of $A_n(M,LM^{-4n})$
is quadratic and exactly agrees with $\NP'_n$ for all $n>1$.

The above discussion implies that $P_n(M,L)$ is a rational mutiple of 
$A_n(M,LM^{-4n})$. Since their leading coefficients (with respect to $L$) 
agree, they are equal. This proves Theorem \ref{thm.1} for $n>1$. The
case of $n<-1$ is similar.
\qed

\subsection{Proof of Theorem \ref{thm.2}}
\lbl{sub.thm2}

Using Equations \eqref{eq.Pn} and \eqref{eq.Qn}, convert 
Equation \eqref{eq.duality} into

\begin{equation}
\lbl{eq.Pduality}
Q_{-n}(\sqrt{M},L/M^{5}) 
= (-L)^n M^{n+13} P_n(i \sqrt{L}, L^{5}/M).
\end{equation}
Note that, under the substitution $(M,L) \mapsto (i/\sqrt{L}, L^{2n+5}/M$), 
$LM^{4n}$ becomes $L^5/M$. Similarly, $LM^{-4n}$ becomes $L/M^5$ under the 
substitution $(M,L) \mapsto (\sqrt{M}, LM^{2n-5})$. 

It is straightforward to verify equation \eqref{eq.Pduality} for 
$n = 2, 3, 4, 5$. For $n \geq 6$, we use induction.
Let $c_k^-$ denote the result of applying the 
substitutions $(M,L) \mapsto (\sqrt{M}, L/M^5)$ to the $c_k$ 
coefficients in the recursions \eqref{eq.recP} and \eqref{eq.recQ}. For 
example,
$$c_0^- = \frac{L^4 (1 + L)^4 (1 - M)^4}{M^2}.$$
Similarly, define $c_k^+$ to be the result of the 
substitution $(M,L) \mapsto (i /\sqrt{L}, L^5/M)$ to $c_k$.
It is easy to verify that for $k = 0,1,2,3$,
$$
\frac{c_k^-}{c_4^-}  (-LM)^{k-4} = \frac{c_k^+}{c_4^+}.
$$ Then,
\begin{eqnarray*}
Q_{-n}(\sqrt{M},L/M^{5}) & = &
 - \frac{1}{c_4^-} \sum_{k=0}^3 c_k^- Q_{-n+4-k} (\sqrt{M}, L/M^5) \\
& = &  - \frac{1}{c_4^-} \sum_{k=0}^3 c_k^-  (-L)^{n-4+k} M^{n-4+k+13} P_{n-4+k} 
(i \sqrt{L}, L^{5}/M) \\
& = &  -  (-L)^n M^{n+13} \sum_{k=0}^3 \frac{c_k^-}{c_4^-}  (-LM)^{k-4} P_{n-4+k} 
(i \sqrt{L}, L^{5}/M) \\
& = & -  (-L)^n M^{n+13} \sum_{k=0}^3 \frac{c_k^+}{c_4^+} P_{n-4+k} 
(i \sqrt{L}, L^{5}/M) \\
& = & (-L)^n M^{n+13} P_{n}(i \sqrt{L}, L^{5}/M).
\end{eqnarray*}

By induction, equation \eqref{eq.Pduality} holds for all $n > 1$
proving Theorem~\ref{thm.2}.
\qed


\appendix

\section{The coefficients $c_k$ and the initial conditions for 
$P_n$ and $Q_n$}
\lbl{sec.ck}


{\small
\begin{eqnarray*}
c_4 &=&
M^4
\\ 
c_3 &=&
1+M^4+2 L M^{12}+L M^{14}-L M^{16}+L^2 M^{20}-L^2 M^{22}-2 L^2 M^{24}-L^3 M^{32}-L^3 M^{36}
\\
c_2 &=& 
\left(-1+L M^{12}\right) \left(-1-2 L M^{10}-3 L M^{12}+2 L M^{14}-L^2 M^{16}+2 L^2 M^{18}-4 L^2 M^{20}-2 L^2 M^{22}+3 L^2 M^{24}
\right. \\
& & \left.
-3 L^3 M^{28}+2 L^3 M^{30}+4 L^3 M^{32}-2 L^3 M^{34}+L^3 M^{36}-2 L^4 M^{38}+3 L^4 M^{40}+2 L^4 M^{42}+L^5 M^{52}\right)
\\
c_1 &=&
-L^2 (-1+M)^2 M^{16} (1+M)^2 \left(1+L M^{10}\right)^2 \left(-1-M^4-2 L M^{12}-L M^{14}+L M^{16}-L^2 M^{20}+L^2 M^{22}
\right. \\
& & \left.
+2 L^2 M^{24}+L^3 M^{32}+L^3 M^{36}\right)
\\
c_0 &=& 
L^4 (-1+M)^4 M^{36} (1+M)^4 \left(1+L M^{10}\right)^4
\end{eqnarray*}
}

{\small
\begin{eqnarray*}
P_0 &=&
\frac{\left(-1+L M^{12}\right) \left(1+L M^{12}\right)^2}{\left(1+L M^{10}\right)^3}
\\ 
P_1 &=&
\frac{\left(-1+L M^{11}\right)^2 \left(1+L M^{11}\right)^2}{1+L M^{10}}
\\ 
P_2 &=&
-1+L M^8-2 L M^{10}+L M^{12}+2 L^2 M^{20}+L^2 M^{22}-L^4 M^{40}-2 L^4 M^{42}-L^5 M^{50}+2 L^5 M^{52}-L^5 M^{54}
\\ & & 
+L^6 M^{62}
\\
P_3 &=&
\left(-1+L M^{12}\right) \left(-1+L M^4-L M^6+2 L M^8-5 L M^{10}+L M^{12}+5 L^2 M^{16}-4 L^2 M^{18}+L^2 M^{22}+L^3 M^{26}
\right. \\
& & \left.
+3 L^3 M^{30}+2 L^3 M^{32}-2 L^4 M^{36}-3 L^4 M^{38}+3 L^4 M^{40}+2 L^4 M^{42}-2 L^5 M^{46}-3 L^5 M^{48}-L^5 M^{52}
\right. \\
& & \left.
-L^6 M^{56}+4 L^6 M^{60}-5 L^6 M^{62}-L^7 M^{66}+5 L^7 M^{68}-2 L^7 M^{70}+L^7 M^{72}-L^7 M^{74}+L^8 M^{78}\right)
\end{eqnarray*}
}

{\small
\begin{eqnarray*}
Q_{0} &=& 
-\frac{\left(-1+L M^{12}\right) \left(1+L M^{12}\right)^2}{L^3 (-1+M)^3 M^4 (1+M)^3}\\
Q_{-1} &=&
-\frac{M^{12} \left(1+L M^{14}\right)^2}{L (-1+M) (1+M)}
\\
Q_{-2} &=&
M^{20} \left(1-L M^8+2 L M^{10}+2 L M^{12}-L M^{16}+L M^{18}+L^2 M^{20}-L^2 M^{22}+2 L^2 M^{26}+2 L^2 M^{28}
\right. \\
& & \left.
-L^2 M^{30}+L^3 M^{38}\right)
\\
Q_{-3} &=&
M^{16} \left(-1+L M^{12}\right) \left(1+L M^{10}+5 L M^{12}-L M^{14}-2 L M^{16}+2 L M^{18}-L M^{20}+2 L^2 M^{20}+L M^{22}
\right. \\
& & \left.
+4 L^2 M^{22}+3 L^2 M^{26}-3 L^2 M^{28}-L^3 M^{28}+5 L^3 M^{30}+5 L^2 M^{32}-L^2 M^{34}-3 L^3 M^{34}+3 L^3 M^{36}
\right. \\
& & \left.
+4 L^3 M^{40}+L^4 M^{40}+2 L^3 M^{42}-L^4 M^{42}+2 L^4 M^{44}-2 L^4 M^{46}-L^4 M^{48}+5 L^4 M^{50}+L^4 M^{52}+L^5 M^{62}\right)
\end{eqnarray*}
}

\section{The coefficients $\ga_k$ and the initial conditions for $R_n$}
\lbl{sec.Rn}

{\small
\begin{eqnarray*}
\ga_4 &=&
L^4 (-1+M)^4 M^{36} (1+M)^4 \left(1+L M^{10}\right)^4
\\
\ga_3 &=&
L^3 (-1+M)^3 M^{24} (1+M)^3 \left(1+L M^{10}\right)^3 \left(-1-M^4-2 L M^{12}-L M^{14}+L M^{16}-L^2 M^{20}+L^2 M^{22}
\right. \\
& & \left.
+2 L^2 M^{24}+L^3 M^{32}+L^3 M^{36}\right)
\\
\ga_2 &=&
L^2 (-1+M)^2 M^{16} (1+M)^2 \left(1+L M^{10}\right)^2 \left(-1+L M^{12}\right) \left(-1-2 L M^{10}-3 L M^{12}+2 L M^{14}-L^2 M^{16}
\right. \\
& & \left.
+2 L^2 M^{18}-4 L^2 M^{20}-2 L^2 M^{22}+3 L^2 M^{24}-3 L^3 M^{28}+2 L^3 M^{30}+4 L^3 M^{32}-2 L^3 M^{34}+L^3 M^{36}
\right. \\
& & \left.
-2 L^4 M^{38}+3 L^4 M^{40}+2 L^4 M^{42}+L^5 M^{52}\right)
\\
\ga_1 &=& \ga_3
\\
\ga_0 &=& \ga_4
\end{eqnarray*}
}
Let $P_n$ for $n=0,\dots,3$ be as in Appendix \ref{sec.ck}. Then,
\begin{equation}
\lbl{eq.Rinit}
R_n=P_n b^n
\end{equation}
for $n=0,\dots,3$ where $b$ is given by Equation \eqref{eq.b}.

\bibliographystyle{hamsalpha}\bibliography{biblio}
\end{document}